\documentclass[a4paper,12pt]{article}
\usepackage{comment}
\usepackage{cite}
\usepackage{amsmath}
\usepackage{amssymb}
\usepackage{amsfonts}
\usepackage[T1]{fontenc}
\usepackage[utf8]{inputenc}
\usepackage{graphicx}
\usepackage{fancyhdr}
\usepackage{float}
\usepackage{xcolor}
\usepackage{authblk}
\usepackage{mathrsfs}
\usepackage{empheq}
\usepackage[hyphens]{url}
\usepackage{hyperref} 
\usepackage[]{breakurl}

\usepackage{geometry}
\begin{document}

\title{Kullback-Leibler divergence for the Fr\'echet extreme-value distribution}

\author[$\dagger$]{Jean-Christophe {\sc Pain}$^{1,2,}$\footnote{jean-christophe.pain@cea.fr}\\
\small
$^1$CEA, DAM, DIF, F-91297 Arpajon, France\\
$^2$Universit\'e Paris-Saclay, CEA, Laboratoire Mati\`ere en Conditions Extr\^emes,\\ 
91680 Bruy\`eres-le-Ch\^atel, France
}

\maketitle

\begin{abstract}
We derive a closed-form solution for the Kullback-Leibler divergence between two Fr\'echet extreme-value distributions. The resulting expression is rather simple and involves the Euler-Mascheroni constant.
\end{abstract}

\section{Introduction}

In probability theory and statistics, the Fr\'echet law is a particular case of a generalized extremum law (which is useful to represent phenomena with extreme values) on the same footing as the Gumbel \cite{Gumbel1958} and Weibull \cite{Papoulis2002} laws. The law was named after Maurice Fr\'echet \cite{Frechet1927}. 

In hydrology for instance, the Fr\'echet distribution is applied to extreme events such as annually maximum one-day rainfalls and river discharges \cite{Coles2001}.

The one-parameter Fr\'echet distribution (probability distribution function) reads
\begin{equation}
f(x;\alpha)=\alpha~x^{-\alpha-1}e^{-x^{-\alpha}}, 
\end{equation}
with $\alpha>0$. Its moments are respectively (see for instance Ref. \cite{Muraleedharan2009}):
\begin{equation}
\mu_{k}=\int_{0}^{\infty}x^{k}f(x;\alpha)\,dx=\int_{0}^{\infty }t^{-{\frac {k}{\alpha }}}e^{-t}\,dt=\Gamma \left(1-{\frac {k}{\alpha}}\right),
\end{equation}
where $\Gamma$ is the usual Gamma function and $k\geq 1$. The moments $\mu_k$ are defined for $k<\alpha$. The centered moments are defined as
\begin{equation}
\mu_{k,c}=\int_{0}^{\infty}(x-\mu_1)^{k}f(x;\alpha)\,dx
\end{equation}
and the reduced centered moments
\begin{equation}
\zeta_{k,c}=\frac{\mu_{k,c}}{(\mu_{2,c})^{n/2}}.
\end{equation}
The skewness of the Fr\'echet distribution is
\begin{equation}
\zeta_{3,c}=\frac{\Gamma \left(1-{\frac {3}{\alpha }}\right)-3\Gamma \left(1-{\frac {2}{\alpha }}\right)\Gamma \left(1-{\frac {1}{\alpha }}\right)+2\Gamma ^{3}\left(1-{\frac {1}{\alpha }}\right)}{\sqrt {\left(\Gamma \left(1-{\frac {2}{\alpha }}\right)-\Gamma ^{2}\left(1-{\frac {1}{\alpha }}\right)\right)^{3}}}
\end{equation}
for $\alpha>3$ and $+\infty$ otherwise. The excess kurtosis (kurtosis minus three) reads
\begin{equation}
\zeta_{4,c}-3=-6+{\frac {\Gamma \left(1-{\frac {4}{\alpha }}\right)-4\Gamma \left(1-{\frac {3}{\alpha }}\right)\Gamma \left(1-{\frac {1}{\alpha }}\right)+3\Gamma ^{2}\left(1-{\frac {2}{\alpha }}\right)}{\left[\Gamma \left(1-{\frac {2}{\alpha }}\right)-\Gamma ^{2}\left(1-{\frac {1}{\alpha }}\right)\right]^{2}}}
\end{equation}
for $\alpha >4$ and $\infty$ otherwise. 

The Kullback-Leibler (KL) divergence \cite{Kullback1951,Kullback1959} provides a non-symmetric measure of the similarity of two probability distributions $p$ and $q$. In the case where both distributions are continuous, it is defined as
\begin{equation}
\mathscr{D}_{\mathrm{KL}}(p|q) = \int_{-\infty}^{\infty} p(x) \ln\left(\frac{p(x)}{q(x)}\right)\,dx,
\end{equation}
where $p(x)$ and $q(x)$ denote the corresponding probability densities. The KL divergence is a measure of relative entropy. It can be understood as the information loss if $p$ is modeled by means of $p$. Accordingly, the smaller the KL divergence is, the more similar are $p$ and $q$. 

The KL divergence was already investigated in the framework of Weibull \cite{Yari2013,Bauckhage2013} and Gumbel \cite{Gil2011,Subh2014} distributions.

In the next section, we present the KL divergence for the one-parameter Fr\'echet distribution.

\section{Calculation of the Kullback-Leibler divergence}

Let us consider two Fr\'echet distributions $f(x;\alpha_1)$ and $f(x;\alpha_2)$. The KL divergence reads
\begin{equation}
\mathscr{D}_{\mathrm{KL}}=\int_0^{\infty}f(x;\alpha_1) \ln \left(\frac{f(x;\alpha_1)}{f(x;\alpha_2)}\right)\,dx. 
\end{equation}
We have to calculate
\begin{equation}
\mathscr{D}_{\mathrm{KL}}=\int_{0}^{\infty} f(x;\alpha_1) \ln\left(\frac{f(x;\alpha_1)}{f(x;\alpha_2)}\right)\,dx,
\end{equation}
\emph{i.e.}
\begin{equation}
\mathscr{D}_{\mathrm{KL}}=\int_{0}^{\infty} \alpha_1~x^{-\alpha_1-1}e^{-x^{-\alpha_1}} \ln\left(\frac{\alpha_1~x^{-\alpha_1-1}e^{-x^{-\alpha_1}}}{\alpha_2~x^{-\alpha_2-1}e^{-x^{-\alpha_2}}}\right)\,dx,
\end{equation}
which is equal to
\begin{equation}
\mathscr{D}_{\mathrm{KL}}=\alpha_1\int_{0}^{\infty} x^{-\alpha_1-1}e^{-x^{-\alpha_1}} \ln\left(\frac{\alpha_1}{\alpha_2}~x^{-\alpha_1+\alpha_2}e^{-x^{-\alpha_1}+x^{-\alpha_2}}\right)\,dx,
\end{equation}
yielding
\begin{eqnarray}
\mathscr{D}_{\mathrm{KL}}&=&\alpha_1\ln\left(\frac{\alpha_1}{\alpha_2}\right)\int_{0}^{\infty} x^{-\alpha_1-1}e^{-x^{-\alpha_1}}\,dx\nonumber\\
& &+\alpha_1(\alpha_2-\alpha_1)\int_{0}^{\infty} x^{-\alpha_1-1}e^{-x^{-\alpha_1}}\ln(x)\,dx\nonumber\\
& &-\alpha_1\int_{0}^{\infty} x^{-\alpha_1-1}e^{-x^{-\alpha_1}}(x^{-\alpha_1}+x^{-\alpha_2})\,dx
\end{eqnarray}
or
\begin{eqnarray}\label{sum}
\mathscr{D}_{\mathrm{KL}}&=&\alpha_1\ln\left(\frac{\alpha_1}{\alpha_2}\right)\int_{0}^{\infty} x^{-\alpha_1-1}e^{-x^{-\alpha_1}} \,dx\nonumber\\
& &+\alpha_1(\alpha_2-\alpha_1)\int_{0}^{\infty} x^{-\alpha_1-1}e^{-x^{-\alpha_1}}\ln(x)\,dx\nonumber\\
& &-\alpha_1\int_{0}^{\infty} x^{-2\alpha_1-1}e^{-x^{-\alpha_1}} \,dx\nonumber\\
& &+\alpha_1\int_{0}^{\infty} x^{-\alpha_2-\alpha_1-1}e^{-x^{-\alpha_1}} \,dx.
\end{eqnarray}
Let us make the change of variable $t=e^{-x^{-\alpha_1}}$. The first integral in Eq. (\ref{sum}) becomes
\begin{equation}\label{first}
\int_{0}^{\infty} x^{-\alpha_1-1}e^{-x^{-\alpha_1}} \,dx=\frac{1}{\alpha_1}\int_0^{\infty}e^{-t}\,dt=\frac{1}{\alpha_1}.
\end{equation}
The second integral in Eq. (\ref{sum}) is
\begin{equation}\label{second}
\int_{0}^{\infty} x^{-\alpha_1-1}e^{-x^{-\alpha_1}}\ln(x)\,dx=-\frac{1}{\alpha_1^2}\int_0^{\infty}\ln t~e^{-t}\,dt=\frac{\gamma}{\alpha_1^2},
\end{equation}
where
\begin{equation}
\gamma=-\int_0^{\infty}\ln t~e^{-t}\,dt=-\Gamma'(1)
\end{equation}
is the usual Euler-Mascheroni constant.

The fourth integral in Eq. (\ref{sum}) is
\begin{equation}\label{fourth}
\int_{0}^{\infty} x^{-\alpha_1-\alpha_2-1}e^{-x^{-\alpha_1}} \,dx=\frac{1}{\alpha_1}\int_0^{\infty}t^{\alpha_2/\alpha_1}~e^{-t}\,dt=\frac{1}{\alpha_1}\Gamma\left(\frac{\alpha_1+\alpha_2}{\alpha_1}\right).
\end{equation}
The third integral in Eq. (\ref{sum}) is given by expression (\ref{fourth}) setting $\alpha_2=\alpha_1$. It reads therefore
\begin{equation}\label{third}
\int_{0}^{\infty} x^{-2\alpha_1-1}e^{-x^{-\alpha_1}} \,dx=\frac{1}{\alpha_1}\int_0^{\infty}t~e^{-t}\,dt=\frac{1}{\alpha_1}\Gamma\left(2\right)=\frac{1}{\alpha_1}.
\end{equation}
Inserting expressions (\ref{first}), (\ref{second}), (\ref{fourth}) and (\ref{third}) in Eq. (\ref{sum}) gives the final result
\begin{empheq}[box=\fbox]{align}
\mathscr{D}_{\mathrm{KL}}=\ln\left(\frac{\alpha_1}{\alpha_2}\right)+\frac{(\alpha_2-\alpha_1)}{\alpha_1}~\gamma-1+\Gamma\left(\frac{\alpha_1+\alpha_2}{2}\right).
\end{empheq}

\section{Conclusion}

In this note, we derived a closed-form solution for the Kullback-Leibler divergence between two Fr\'echet extreme-value distributions.
The present results can be easily extended to the generalized Fr\'echet law, obtained by introducing a parameter $m$ (position of the maximum) and a scale parameter $s>0$. The cumulative distribution function reads then $e^{-\left(\frac{x-m}{s}\right)^{-\alpha}}$ (instead of $e^{-x^{\alpha}}$ in the one-parameter case) if $x>m$ and 0 otherwise. It corresponds to the probability distribution function
\begin{equation}
g(x;\alpha,s,m)=\frac{\alpha}{s}\left(\frac{x-m}{s}\right)^{-1-\alpha}e^{-(\frac{x-m}{s})^{-\alpha}}.
\end{equation}


\begin{thebibliography}{99}

\bibitem{Gumbel1958} 
E. J. Gumbel, {\it Statistics of extremes} (New York: Columbia University Press, 1958).

\bibitem{Papoulis2002} 
A. Papoulis and S. U. Pillai, {\it Probability, random variables, and stochastic processes} (Boston: McGraw-Hill, 2002).

\bibitem{Frechet1927} 
M. Fr\'echet, {\it Sur la loi de probabilit\'e de l'\'ecart maximum}, Ann. Soc. Polon. Math. {\bf 6}, 3 (1927) [in french].

\bibitem{Coles2001}
S. Coles, {\it An introduction to statistical modeling of extreme values} (Springer-Verlag, 2001). 

\bibitem{Muraleedharan2009} G. Muraleedharan, C. Guedes Soares and C. Lucas, {\it Characteristic and moment generating functions of generalized extreme value distribution (GEV)}, chapter 13, in: {\it Sea level rise, coastal engineering, shorelines and tides}, edited by Linda L. Wright (Nova Science Publishers, Inc, 2009).

\bibitem{Kullback1951} 
S. Kullback et R. Leibler, {\it On information and sufficiency}, Ann. Math. Stat. {\bf 22}, 79 (1951).

\bibitem{Kullback1959} 
S. Kullback, {\it Information theory and statistics} (New York, John Wiley and Sons, 1959).

\bibitem{Yari2013}
G. Yari, A. Mirhabibi and A. Saghafi, {\it Estimation of the Weibull parameters by Kullback-Leibler divergence of survival functions}, Appl. Math. Inf. Sci. {\bf 7}, 187 (2013).

\bibitem{Bauckhage2013}
Ch. Bauckhage, {\it Computing the Kullback-Leibler divergence between two Weibull distributions}, arXiv 1310.3713 (2013).\\
\url{https://arxiv.org/pdf/1310.3713.pdf}

\bibitem{Gil2011}
M. Gil, {\it On R\'enyi divergence measures for continuous alphabet sources}, Masters of Science Thesis, Department of Mathematics and Statistics, Queen's University, Kingston, Ontario, Canada (2011).

\bibitem{Subh2014}
S. A. Al-Subh, {\it Goodness of fit test for Gumbel distribution based on Kullback-Leibler information using several different estimators}, Appl. Math. Sci. {\bf 8}, 4703 (2014).

\end{thebibliography}
\end{document}